\documentclass[14pt, draft]{article}
\usepackage[cp1251]{inputenc}
\usepackage[ukrainian, russian, english]{babel}

\usepackage{ amsfonts, amssymb}
\usepackage{amsmath}
\usepackage{amsthm}

\begin{document}
\selectlanguage{ukrainian} \thispagestyle{empty}
 \pagestyle{myheadings}              

УДК 517.5 \vskip 3mm

\noindent \bf А.С. Сердюк  \rm (Інститут математики НАН України, Київ) \\
\noindent \bf У.З. Грабова  \rm (Східноєвропейський національний університет імені Лесі Українки, Луцьк)

\noindent\bf A.S. Serdyuk   \rm (Institute of Mathematics of The National Academy of Sciences
of Ukraine, Kiev)  \\
 \noindent\bf U.Z. Grabova \rm (Lesja Ukrainka East European National University, Lutsk)

 \vskip 5mm

\centerline \textbf{ Порядкові оцінки найкращих наближень і
наближень сумами Фур'є  класів $(\psi,\beta)$ -- диференційовних
функцій}

\vskip 5mm

\centerline \textbf{Order estimation  of the best approximations
and of the approximations by Fourier sums of classes of
$(\psi,\beta)$--differentiable functions }

\vskip 5mm

 \rm Встановлено точні за порядком оцінки найкращих рівномірних
наближень  тригонометричними поліномами  на класах
$C^{\psi}_{\beta,p}$  ---
 $2\pi$--періодичних неперервних функцій $f$, які зображуються згортками функцій, що належать одиничним кулям просторів $L_{p}$, $1\leq p<\infty$,
  з фіксованими твірними ядрами $\Psi_{\beta}\subset L_{p'}$,
  $\frac{1}{p}+\frac{1}{p'}=1$, коефіцієнти Фур'є яких спадають до
  нуля приблизно як степеневі функції. Точні порядкові оцінки
  найкращих наближень встановлено також і в $L_{p}$--метриці,
  $1<p\leq\infty$,  для класів $L^{\psi}_{\beta,1}$
---  $2\pi$--періодичних функцій $f$ еквівалентних відносно міри Лебега до згорток ядер $\Psi_{\beta}\subset L_{p}$ із функціями з одиничної кулі простору
$L_{1}$.
 Показано, що в усіх розглядуваних випадках порядки  найкращих наближень  реалізують суми
 Фур'є.

\vskip 5mm

 \rm There were established the exact--order estimations of the best uniform approximations by
 the trigono\-metrical polynoms  on the $C^{\psi}_{\beta,p}$
 classes of $2\pi$--periodic continuous functions $f$, which are
 defined by the convolutions  of the functions, which belong to
 the unit ball in $L_{p}$, $1\leq p<\infty$ spaces with
 generating fixed kernels $\Psi_{\beta}\subset L_{p'}$,
 $\frac{1}{p}+\frac{1}{p'}=1$, whose Fourier coefficients
 decreasing to zero approximately as power
  functions. The exact order estimations were also established in
  $L_{p}$--metrics, $1<p\leq\infty$ for $L^{\psi}_{\beta,1}$
  classes of $2\pi$--periodic functions $f$, which are equivalent
  by means of Lebesque measure to the convolutions of $\Psi_{\beta}\subset
  L_{p}$ kernels with the functions that belong to the unit ball in
  $L_{1}$ space. We showed that in investigating cases the orders of best
approximations are realized by Fourier sums.

\newpage

Нехай $L_{p}$, $1\leq p\leq\infty$, --- простір
$2\pi$--періодичних сумовних  функцій $f$ зі скінченною нормою $\|
f\|_{p}$, де при $p\in [1,\infty)$ \
$\|f\|_{p}=\Big(\int\limits_{0}^{2\pi}|f(t)|^{p}dt\Big)^{\frac{1}{p}}$,
а при $p=\infty$ \
$\|f\|_{\infty}=\mathop{\rm{ess}\sup}\limits_{t}|f(t)|$,
 $C$ --- простір $2\pi$--періодичних неперервних функцій, у
якому норма задається рівністю $\|f\|_{C}=\max\limits_{t}|f(t)|$.

Нехай, далі  $L^{\psi}_{\beta,p}$ --- клас $2\pi$--періодичних
функцій $f(x)$, котрі майже для всіх $x\in\mathbb{R}$
представляються згортками

\begin{equation}\label{zgo}
f(x)=\frac{a_{0}}{2}+\frac{1}{\pi}\int\limits_{-\pi}^{\pi}\Psi_{\beta}(x-t)\varphi(t)dt,
\ a_{0}\in\mathbb{R},  \ \varphi\perp 1,
\end{equation}
де  $\|\varphi\|_{p}\leq 1$, $1\leq p\leq\infty$,
$\Psi_{\beta}(t)$
--- сумовна на $(0,2\pi)$ функція, ряд Фур'є якої має вигляд
\begin{equation}\label{rff}
\sum\limits_{k=1}^{\infty}\psi(k)\cos
\big(kt-\frac{\beta\pi}{2}\big), \ \psi(k)>0, \  \beta\in
    \mathbb{R}.
\end{equation}
 Функцію
$\varphi$ в зображенні (\ref{zgo}), згідно з О.І. Степанцем
[\ref{S1}, с. 132 ], називають $(\psi,\beta)$--похідною функції
$f$ і позначають через $f^{\psi}_{\beta}$.

При $\psi(k)=k^{-r}$ класи $L^{\psi}_{\beta,p}$ перетворюються у
відомі класи  Вейля--Надя $W^{r}_{\beta,p}$, а їх
$(\psi,\beta)$--похідні $f^{\psi}_{\beta}$ майже скрізь
співпадають з похідними в сенсі  Вейля--Надя $f^{r}_{\beta}$,
останні при $r=\beta$, \ $r\in\mathbb{N}$ майже скрізь збігаються
зі звичайними $r$--ми похідними функції $f$.

Якщо твірне ядро $\Psi_{\beta}$ класу $L^{\psi}_{\beta,p}$
задовольняє включенню $\Psi_{\beta}\in L_{p'}$,
$\frac{1}{p}+\frac{1}{p'}=~1$, то $L^{\psi}_{\beta,p}\subset
L_{\infty}$, $1\leq p\leq\infty$, а згортки виду (\ref{zgo}) є
неперервними функціями \linebreak(див. твердження 3.8.1 роботи
[\ref{S1}, с. 137]). Тому клас усіх функцій $f$ виду (\ref{zgo}),
для яких $\|\varphi\|_{p}\leq 1$, $\Psi_{\beta}\in L_{p'}$  будемо
позначати через $C^{\psi}_{\beta,p}$.

У випадку, якщо $\Psi_{\beta}\in L_{p}$, $1\leq p\leq\infty$,  то
(див., наприклад, [\ref{Korn}, с. 71]) має місце включення
$L^{\psi}_{\beta,1}\subset L_{p}$.

В даній роботі розглядається задача про знаходження точних
порядкових оцінок  функціональних класів $L^{\psi}_{\beta,p}$
 та $L^{\psi}_{\beta,1}$ сумами Фур'є
$S_{n-1}(t)$ порядку $n-1$ у метриках $L_{\infty}$ та $L_{p}$
відповідно

\begin{equation}\label{och1}
{\cal E}_{n}(L^{\psi}_{\beta,p})_{L_{\infty}}=\sup\limits_{f\in
L^{\psi}_{\beta,p}}\|f(\cdot)-S_{n-1}(f;\cdot)\|_{\infty}, \ \
1\leq p<\infty,
\end{equation}
\begin{equation}\label{och11}
{\cal E}_{n}(L^{\psi}_{\beta,1})_{L_{p}}=\sup\limits_{f\in
L^{\psi}_{\beta,1}}\|f(\cdot)-S_{n-1}(f;\cdot)\|_{p}, \ \ \ 1<
p\leq\infty,
\end{equation}
а також задача про знаходження точних порядкових оцінок найкращих
наближень класів $L^{\psi}_{\beta,p}$ та $L^{\psi}_{\beta,1}$ в
метриках $L_{\infty}$ та $L_{p}$ відповідно
\begin{equation}\label{och2}
{E}_{n}(L^{\psi}_{\beta,p})_{L_{\infty}}=\sup\limits_{f\in
L^{\psi}_{\beta,p}}\inf\limits_{t_{n-1}\in\mathcal{T}_{2n-1}}\|f(\cdot)-t_{n-1}(\cdot)\|_{\infty},
\ \ 1\leq p<\infty,
\end{equation}
\begin{equation}\label{och22}
{E}_{n}(L^{\psi}_{\beta,1})_{L_{p}}=\sup\limits_{f\in
L^{\psi}_{\beta,1}}\inf\limits_{t_{n-1}\in\mathcal{T}_{2n-1}}\|f(\cdot)-t_{n-1}(\cdot)\|_{p},
\ \ \ 1< p\leq\infty,
\end{equation}
де $\mathcal{T}_{2n-1}$ --- підпростір усіх тригонометричних
поліномів $t_{n-1}$ порядку не вищого за $n-1$.

Зрозуміло, що у випадку класів $C^{\psi}_{\beta,p}$ норму
$\|\cdot\|_{\infty}$ в (\ref{och1}) і (\ref{och2}) слід замінити
на $\|\cdot\|_{C}$
 і при цьому  ${\cal E}_{n}(C^{\psi}_{\beta,p})_{C}={\cal E}_{n}(L^{\psi}_{\beta,p})_{L_{\infty}}$, ${E}_{n}(C^{\psi}_{\beta,p})_{C}={E}_{n}(L^{\psi}_{\beta,p})_{L_{\infty}}$.

Для класів Вейля--Надя   $W^{r}_{\beta,p}$, $\beta\in\mathbb{R}$,
$1\leq p\leq\infty$ точні порядкові оцінки величин
\linebreak(\ref{och1}) -- (\ref{och22}) відомі і мають вигляд
(див., наприклад, [\ref{T}, с. 47--49])
\begin{equation}\label{teml1}
{\cal E}_{n}(W^{r}_{\beta,p})_{\infty}\asymp n^{-r+\frac{1}{p}}, \
\ \ \  1\leq p<\infty, \ \ \ \ r>\frac{1}{p}, \ \ \ \ \ \ \ \ \ \
\ \ \ \ \ \ \ \ \ \
\end{equation}
\begin{equation}\label{teml3}
{E}_{n}(W^{r}_{\beta,p})_{\infty}\asymp n^{-r+\frac{1}{p}},  \ \ \
\   1\leq p\leq\infty, \ \ \ \  r>\frac{1}{p},\ \ \ \ \ \ \ \ \ \
\ \ \ \ \ \ \ \ \ \
\end{equation}
\begin{equation}\label{teml14}
{\cal E}_{n}(W^{r}_{\beta,1})_{p}\asymp n^{-r+\frac{1}{p'}}, \ \ \
 \ \
 1< p\leq\infty, \ \ \ \ r>\frac{1}{p'}, \ \ \ \
\frac{1}{p}+\frac{1}{p'}=1,
\end{equation}
\begin{equation}\label{teml15}
{E}_{n}(W^{r}_{\beta,1})_{p}\asymp n^{-r+\frac{1}{p'}}, \ \ \ \
1\leq p\leq\infty, \ \ \ \ r>\frac{1}{p'}, \ \ \ \
\frac{1}{p}+\frac{1}{p'}=1,
\end{equation}
\begin{equation}\label{teml2}
{\cal E}_{n}(W^{r}_{\beta,p})_{p}\asymp n^{-r}\ln n,  \ \ \ \ \ \
  p=1,\infty, \ \ \ \ r>0, \ \ \ \ n\in \mathbb{N}\setminus \{1\}.
\end{equation}
У формулах (\ref{teml1}) -- (\ref{teml2}) і надалі під записом
$A_{n}\asymp B_{n}$ будемо розуміти існування додатних сталих
$K_{1}$ і $K_{2}$ таких, що $K_{1}B_{n}\leq A_{n}\leq K_{2}B_{n}$,
$n\in\mathbb{N}$.

Щодо випадку $p=1, \infty$ зауважимо, що завдяки роботам А.М.
Колмогорова  [\ref{Kol}], В.Т.~ Пінкевича [\ref{Pin}],  С.М.
Нікольського [\ref{Nik}], [\ref{Nik2}],  А.В. Єфімова [\ref{Ef}]
та С.О. Теляковського~ [\ref{Tel}] для величин ${\cal
E}_{n}(W^{r}_{\beta,\infty})_{\infty}$ та ${\cal
E}_{n}(W^{r}_{\beta,1})_{1} $ при $r>0$, $\beta\in\mathbb{R}$
відомі асимптотичні рівності при $n\rightarrow\infty$.

Що ж стосується найкращих наближень
${E}_{n}(W^{r}_{\beta,\infty})_{\infty}$ та
${E}_{n}(W^{r}_{\beta,1})_{1}$, то завдяки роботам Ж. Фавара
[\ref{Fav},  \ref{Fav1}], В.К. Дзядика [\ref{Dz}], [\ref{Dz1}],
С.Б. Стєчкіна [\ref{ST}] та Сунь~ Юн--шена [\ref{SU}] встановлені
точні значення цих величин при усіх $n\in\mathbb{N}$, $r>0$ і
$\beta\in \mathbb{R}$. Точні значення величин ${\cal
E}_{n}(W^{r}_{\beta,p})_{\infty}$ відомі також у випадку $p=2$
[\ref{Bab}].

На класах $L^{\psi}_{\beta,p}$ точні порядкові оцінки величин ${\cal
E}_{n}(L^{\psi}_{\beta,p})_{s}$ та
$E_{n}(L^{\psi}_{\beta,p})_{s}$ у випадку, коли $\psi(k)k^{\frac{1}{p}-\frac{1}{s}}$
монотонно незростають і $\frac{\psi(k)}{\psi(2k)}\leq K<\infty$, $k\in
\mathbb{N}$ були знайдені у
роботі О.І.~Степанця та О.К. Кушпеля [\ref{Kuchpel}] при довільних
$1< p$, $s<\infty$.

 Крім того, точні порядкові оцінки величин ${\cal
E}_{n}(L^{\psi}_{\beta,p})_{s}$ та $E_{n}(L^{\psi}_{\beta,p})_{s}$ одержані О.І.~ Степанцем [\ref{Step monog 1987}, с. 48] при довільних
$1\leq p$, $s\leq\infty$ за умови $\psi\in {\mathfrak
M^{'}_{\infty}}$ (в цьому випадку $\psi(k)$ спадають до нуля не
повільніше ніж члени деякої геометричної прогресії), а також при
довільних $1<p$, $s<\infty$ за умови $\psi\in {\mathfrak
M^{''}_{\infty}}$ (в зазначеному випадку $\psi(k)$ спадають  до нуля швидше довільної
степеневої функції, але не швидше за деяку геометричну прогресію).
Згодом В.С.~Романюк [\ref{Rom}] у випадку $\psi\in {\mathfrak
M^{''}_{\infty}}$ доповнив згадані результати О.І. Степанця,
встановивши точні порядки величин ${\cal
E}_{n}(C^{\psi}_{\beta,p})_{C}$, $1<p<\infty$ (для величин
$E_{n}(C^{\psi}_{\beta,p})_{C}$ при $\psi\in {\mathfrak
M^{''}_{\infty}}$ питання про точні порядкові оцінки до цих пір
залишається відкритим). Зазначимо також, що при $p=2$ точні
значення величин ${\cal E}_{n}(C^{\psi}_{\beta,p})_{C}$ для всіх
$n\in\mathbb{N}$, $\beta\in \mathbb{R}$
 за умови збіжності ряду
$\sum\limits_{k=1}^{\infty}\psi^{2}(k)$ знайдені у роботі А.С.
Сердюка та І.В. Соколенка [\ref{Serduk}]. Задача про точні
значення величин ${\cal E}_{n}(L^{\psi}_{\beta,2})_{2}$ та
$E_{n}(L^{\psi}_{\beta,2})_{2}$ повністю розв'язана у роботі О.І.~Степанця та О.К.~Кушпеля~[\ref{Kuchpel}].

При $p=1,\infty$ результати, що містять асимптотично точні оцінки
величин ${\cal E}_{n}(L^{\psi}_{\beta,p})_{p}$, а також точні
порядкові оцінки величин $E_{n}(L^{\psi}_{\beta,p})_{p}$ в
залежності від швидкості прямування до нуля послідовності
$\psi(k)$ при $k\rightarrow \infty$ найбільш повно викладені в
монографіях [\ref{S1}], [\ref{S2}].

В даній роботі встановлено точні порядкові оцінки величин
(\ref{och1})--(\ref{och22}) при довільних $\beta \in \mathbb{R}$
у випадку, коли послідовність $\psi(k)$ спадає до нуля не
повільніше і не швидше деяких степеневих функцій. Тим самим
доповнено основні результати роботи [\ref{Kuchpel}] по відшуканню
слабкої асимптотики величин ${\cal E}_{n}(L^{\psi}_{\beta,p})_{s}$
та $E_{n}(L^{\psi}_{\beta,p})_{s}$ на випадки $p=1$ і $s=\infty$.

Перейдемо до точних формулювань.

Вважаючи, що послідовність $\psi(k)$, що визначає клас
$C^{\psi}_{\beta,p}$, є слідом на множині $\mathbb{N }$ деякої
неперервної функції $\psi(t)$ неперервного аргументу $t\geq1$,
позначимо через $\Theta_{p}$, $1\leq p<\infty$, множину монотонно
незростаючих функцій $\psi(t)$, для яких існує стала
$\alpha>\frac{1}{p}$ така, що функція $t^{\alpha}\psi(t)$ майже
спадає, тобто знайдеться додатна стала $K$ така, що
\linebreak$t^{\alpha}_{1}\psi(t_{1})\leq
Kt^{\alpha}_{2}\psi(t_{2})$ для будь--яких $t_{1}>t_{2}\geq 1$;
через $B$ позначимо множину монотонно незростаючих  при $t\geq 1$
додатних функцій $\psi(t)$, для кожної з яких можна вказати
додатну   сталу $K$ таку, що
$$
\frac{\psi(t)}{\psi(2t)}\leq K, \ \  \forall t\geq 1.
$$

Надалі скрізь будемо вважати, що $\psi\in B\cap\Theta_{p}$, $1\leq
p<\infty$. Умова $\psi\in\Theta_{p}, \ 1\leq p<\infty$, як неважко переконатись,
гарантує справедливість включення $\Psi_{\beta}\in L_{p'}$,
$\frac{1}{p}+\frac{1}{p'}=1$ (див., наприклад, [\ref{Bari}, с.
657]). Як випливає \linebreak з [\ref{S1}, с. 165, 175], якщо
$\psi\in B\cap{\mathfrak M}$, де ${\mathfrak M}$ --- множина усіх
опуклих донизу на $[1,\infty)$ функцій $\psi(t)$, таких, що
 $\lim\limits_{t\rightarrow\infty}\psi(t)=0$, то можна вказати таке $r>0$, що при всіх $t\geq 1$ буде виконуватись нерівність $\psi(t)\geq Kt^{-r}$.

  Прикладами функцій $\psi$, що задовольняють умову  $\psi\in B\cap\Theta_{p}$, є, зокрема, функції виду
 $\psi(t)=\frac{1}{t^{r}}$, \ $r>\frac{1}{p}$; \
 $\psi(t)=\frac{1}{t^{r}\ln^{\alpha}(t+c)}$, \ $\alpha\geq0$, \
 $c>0$, \  $r>\frac{1}{p}$, \ $t\geq1$; \ $\psi(t)=\frac{\ln^{\alpha}(t+c)}{t^{r}}$,
 \ $\alpha\geq0$,
 $c>e^{\frac{\alpha}{r}}-1$, \
 $r>\frac{1}{p}$, \ $t\geq1$.

Має місце наступне твердження.

 \noindent \rm {\bf Теорема 1.} \it Нехай  $1<p<\infty$, $\beta\in \mathbb{R}$,
     $\psi\in B\cap
 \Theta_{p}$.
 Тоді  існують додатні величини $K^{(1)}_{\psi,p}$, $K^{(2)}_{\psi,p}$,
 що можуть  залежати лише  від $\psi$ і $p$ такі, що для довільних $n\in \mathbb{N}$
\begin{equation}\label{t11}
K_{\psi,p}^{(2)}\psi(n)n^{\frac{1}{p}}\leq{
E}_{n}\Big(C^{\psi}_{\beta,p}\Big)_{C}\leq{\cal
E}_{n}\Big(C^{\psi}_{\beta,p}\Big)_{C}\leq K_{\psi,p}^{(1)}
\psi(n)n^{\frac{1}{p}}.
\end{equation}

\bf Доведення. \ \rm Для довільної функції $f\in
C^{\psi}_{\beta,p}$, згідно з інтегральним зображенням
(\ref{zgo}), одержимо
\begin{equation}\label{intz}
f(x)-S_{n-1}(f;x)=\frac{1}{\pi}\int\limits_{-\pi}^{\pi}\Psi_{\beta,n}(x-t)\varphi(t)dt,
\end{equation}
де
\begin{equation}\label{fpsib}
\Psi_{\beta,n}(t)=
\sum\limits_{k=n}^{\infty}\psi(k)\cos\big(kt-\frac{\beta\pi}{2}\big).
\end{equation}
Застосовуючи  нерівність Гельдера,  з  рівності  (\ref{intz})
маємо
\begin{equation}\label{ng}
{\cal E}_{n}\big(C^{\psi}_{\beta,p}
\big)_{C}\leq\frac{1}{\pi}\big\|\Psi_{\beta,n}(\cdot)\big\|_{p'}\|
\varphi(\cdot)\|_{p}\leq\frac{1}{\pi}\big\|\Psi_{\beta,n}(\cdot)\big\|_{p'},
\end{equation}
де $\frac{1}{p}+\frac{1}{p'}=1$, $1\leq p<\infty$.

Перетворивши  функцію $\Psi_{\beta,n}(t)$ за допомогою
перетворення Абеля, при довільному $n\in~ \mathbb{N}$, одержимо
\begin{equation}\label{peretA}
\Psi_{\beta,n
}(t)=\sum\limits_{k=n}^{\infty}\Delta\psi(k)D_{k,\beta}(t)-\psi(n)D_{n-1,\beta}(t),
\end{equation}
де $\triangle\psi(k)\mathop{=}\limits^{\rm df}\psi(k)-\psi(k+1)$, а
\begin{equation}\label{dkb}
D_{k,\beta}(t)=\frac{1}{2}\cos\frac{\beta\pi}{2}+\sum\limits_{\nu=1}^{k}\cos\Big(\nu
t-\frac{\beta\pi}{2}\Big)=
$$
$$
=\cos\frac{\beta\pi}{2}\Bigg[\frac{\sin\frac{2k+1}{2}t}{2\sin\frac{t}{2}}\Bigg]+\sin\frac{\beta\pi}{2}\Bigg[\frac{\cos\frac{t}{2}-
\cos\frac{2k+1}{2}t} {2\sin\frac{t}{2}}\Bigg].
\end{equation}
Оскільки
(див., наприклад, [\ref{Z}, с. 13])
\begin{equation}\label{nerdiri}
|D_{k,\beta}(t)|\leq\frac{1}{2}+k, \ \
|D_{k,\beta}(t)|\leq(1+\pi)\Big(\frac{1}{|t|}\Big), \ 0<|t|\leq
\pi,
\end{equation}
то для будь--яких $k\in \mathbb{N}$ і $1<p'<\infty$ маємо
\begin{equation}\label{moddiri}
\int\limits_{-\pi}^{\pi}|D_{k,\beta}(t)|^{p'}dt\leq\int\limits_{0\leq|t|\leq\frac{1}{k}
}(\frac{1}{2}+k)^{p'}dt+\int\limits_{\frac{1}{k}\leq|t|\leq\pi}(1+\pi)^{p'}\frac{dt}{|t|^{p^{'}}}\leq
K_{p'}k^{p'-1},
\end{equation}
де $K_{p'}$ --- стала, що залежить від $p'$. З (\ref{nerdiri}) та
оцінки (\ref{moddiri})
 отримаємо
\begin{equation}\label{ndiri}
\|D_{k,\beta}(t)\|_{p'}\leq K_{p,1}k^{\frac{1}{p}}, \ \
k\in\mathbb{N}, \ \ 1\leq p <\infty, \ \
\frac{1}{p}+\frac{1}{p'}=1, \ \ \beta\in \mathbb{R},
\end{equation}
де $K_{p,1}$ --- стала, що залежить від $p$. Із (\ref{peretA}) та
(\ref{ndiri}) випливає нерівність
\begin{equation}\label{ny}
\|\Psi_{\beta,n}(t)\|_{p'}\leq K_{p,1}
\big(\sum\limits_{k=n}^{\infty}\Delta\psi(k)k^{\frac{1}{p}}+\psi(n)n^{\frac{1}{p}}\big),
\ \ 1\leq p<\infty.
\end{equation}

Для оцінки суми
$\sum\limits_{k=n}^{\infty}\Delta\psi(k)k^{\frac{1}{p}}$ нам буде
корисним наступне твердження.

\noindent \rm {\bf Лема 1.} \it Нехай  $r\in (0,1]$, а
$\psi(k)>0$, монотонно незростає і для неї знайдеться
$\varepsilon>0$ таке, що послідовність $k^{r+\varepsilon}\psi(k)$
майже спадає.  Тоді існує стала $K$, залежна від $\psi$ і $r$
така, що для довільних $n\in \mathbb{N}$
\begin{equation}\label{nr2}
\psi(n)n^{r}\leq\sum\limits_{k=n}^{\infty}\triangle\psi(k)k^{r}\leq
K\psi(n)n^{r}.
\end{equation}

\bf Доведення. \  \rm Оскільки $\psi(k)$ монотонно незростає, то
для $\forall r>0$
\begin{equation}\label{ln1}
\sum\limits_{k=n}^{\infty}\Delta\psi(k)k^{r}\geq
n^{r}\sum\limits_{k=n}^{\infty}\Delta\psi(k)=n^{r}\psi(n).
\end{equation}

Залишається показати, що за виконання умов леми 1 виконується
нерівність
\begin{equation}\label{nA}
\sum\limits_{k=n}^{\infty}\Delta\psi(k)k^{r}\leq
K\psi(n)n^{r}.
\end{equation}

 Застосування  перетворення Абеля дозволяє для будь--яких натуральних
$n\leq M$, і довільного $r\in (0, 1]$ записати рівність
\begin{equation}\label{za}
\sum\limits_{k=n}^{M}\psi(k)k^{r-1}=\sum\limits_{k=n}^{M}\Delta\psi(k)\sum\limits_{\nu=1}^{k}\frac{1}{\nu^{1-r}}-\psi(n)
\sum\limits_{\nu=1}^{n-1}\frac{1}{\nu^{1-r}}+\psi(M+1)\sum\limits_{\nu=1}^{M}\frac{1}{\nu^{1-r}}.
\end{equation}
В силу простих геометричних міркувань неважко переконатись, що для
довільних $m\in \mathbb{N}$ і $r\in (0, 1]$
\begin{equation}\label{sr1}
\frac{1}{r}(m^{r}-1)<\sum\limits_{\nu=1}^{m}\frac{1}{\nu^{1-r}}\leq\frac{1}{r}m^{r}.
\end{equation}
Тому в силу   (\ref{za}) і (\ref{sr1})
\begin{equation}\label{ner1}
\sum\limits_{k=n}^{M}\psi(k)k^{r-1}>\frac{1}{r}\Big(\sum\limits_{k=n}^{M}\Delta\psi(k)(k^{r}-1)-\psi(n)(n-1)^{r}\Big)>
$$
$$
>\frac{1}{r}\Big(\sum\limits_{k=n}^{M}\Delta\psi(k)k^{r}-\psi(n)\big((n-1)^{r}+1\big)\Big).
\end{equation}
При $M\rightarrow\infty$ із (\ref{ner1}) одержуємо
\begin{equation}\label{ner2}
\sum\limits_{k=n}^{\infty}\Delta\psi(k)k^{r}\leq r
\sum\limits_{k=n}^{\infty}\psi(k)k^{r-1}+\psi(n)\big((n-1)^{r}+1\big).
\end{equation}

Оскільки за умовою леми існує $\varepsilon>0$ таке, що
послідовність $k^{r+\varepsilon}\psi(k)$ майже спадає, то
знайдеться стала $K_{1}$ така, що
\begin{equation}\label{fmc}
\psi(k)k^{r+\varepsilon}\leq K_{1}\psi(n)n^{r+\varepsilon}, k=n,
n+1, ...,
\end{equation}
тому
\begin{equation}\label{lozv}
\sum\limits_{k=n}^{\infty}\psi(k)k^{r-1}=\sum\limits_{k=n}^{\infty}\frac{\psi(k)k^{r+\varepsilon}}{k^{1+\varepsilon}}\leq
K_{1}\psi(n)n^{r+\varepsilon}\sum\limits_{k=n}^{\infty}\frac{1}{k^{1+\varepsilon}}\leq
K_{2}\psi(n)n^{r}.
\end{equation}
 Із  (\ref{ner2}) і (\ref{lozv}) випливає (\ref{nr2}). Лему  доведено.

 Оскільки $\psi\in\Theta_{p}$, то, застосувавши лему 1, при $r=\frac{1}{p}$ із (\ref{ng}) та (\ref{ny}) отримуємо нерівність
\begin{equation}\label{ozv}
{ E}_{n}\Big(C^{\psi}_{\beta,p}\Big)_{C}\leq{\cal
E}_{n}\big(C^{\psi}_{\beta,p}\big)_{C}\leq
K_{\psi,p}^{(1)}\psi(n)n^{\frac{1}{p}}, \  1\leq p<\infty,
\end{equation}
$K_{\psi,p}^{(1)}$ --- величина, що залежить від $\psi$ і $p$.

Для того, щоб одержати оцінку знизу розглянемо при заданому $n\in
\mathbb{N}$ функцію
$$
f^{\ast}_{n,\alpha}=
\frac{\alpha\psi(n)}{n^{1-\frac{1}{p}}}\Big(V_{2n}(t)-V_{n}(t)\Big),
\ \alpha>0, \ n\in \mathbb{N},
$$
де  $V_{m}(t)$ --- ядра методу   Валле Пуссена
\begin{equation}\label{uavp}
V_{m}(t)=\frac{1}{m}\sum\limits_{k=m}^{2m-1}D_{k}(t), \ m\in
\mathbb{N},
\end{equation}
$D_{k}(t)$ --- ядра Діріхле
$$
D_{k}(t)=\frac{1}{2}+\sum\limits_{\nu=1}^{k}\cos \nu
t=\frac{\sin\big(k+\frac{1}{2}\big)t}{2\sin\frac{t}{2}}, \ k\in
\mathbb{N}.
$$

Покажемо, спочатку, що при певному виборі значення параметра
$\alpha$ виконується нерівність
\begin{equation}\label{n}
\Big\|\big(f^{\ast}_{n,\alpha}(\cdot)\big)^{\psi}_{\beta}\Big\|_{p}\leq
1, \ 1<p<\infty.
\end{equation}
Для цього скористаємось наступним твердженням роботи [\ref{S2},~
с. 117], в якій встановлено нерівності Бернштейна для $(\psi,
\beta)$-похідних в $L_p$-- метриках для поліномів, тобто нерівності між
$\|(t_{m})^{\psi}_{\beta}\|_{p}$ та $\|t_{m}\|_{p}$,   де $t_{m}$~
--- тригонометричні поліноми порядку $m$.

\noindent \rm {\bf Твердження 1.} \it Нехай $1<p<\infty$, \
$\beta\in  \mathbb{R}$, \ $\psi(k)$ --- довільна незростаюча
послідовність невід'ємних чисел. Тоді для довільного
тригонометричного полінома $t_{m}(\cdot)$ порядку $m$ знайдеться
величина $C_{\psi,p}$, що може залежати тільки   від функції
$\psi(\cdot)$ та числа $p$ така, що
\begin{equation}\label{nerb1}
\big\|\big(t_{m}(\cdot)\big)^{\psi}_{\beta}\big\|_{p}\leq C_{\psi,p}
(\psi(m))^{-1}\|t_{m}(\cdot)\|_{p}.
\end{equation}

\ \rm Оскільки $f^{\ast}_{n,\alpha}$ є тригонометричним поліномом
порядку $4n-1$, то, використовуючи твердження 1, отримаємо
\begin{equation}\label{w1}
\Big\|\big(f^{\ast}_{n,\alpha}(t)\big)^{\psi}_{\beta}\Big\|_{p}\leq\frac{\alpha
C_{\psi,p}}{n^{1-\frac{1}{p}}}\frac{\psi(n)}{\psi(4n-1)}
\|V_{2n}(t)-V_{n}(t)\|_{p}.
\end{equation}

Знайдемо оцінку $\|V_{2n}(t)-V_{n}(t)\|_{p}$. Враховуючи, що
\begin{equation}\label{uaf}
V_{m}(t)=2F_{2m-1}(t)-F_{m-1}(t),
\end{equation}
 де $F_{k}(t)$ --- ядра Фейєра
 $$
 F_{k}(t)=\frac{1}{k+1}\sum\limits_{\nu=0}^{k}D_{\nu}(t)=\frac{1}{k+1}
 \sum\limits_{\nu=0}^{k}\frac{\sin\big(\nu+\frac{1}{2}\big)t}{2\sin\frac{t}{2}}, \ k\in
\mathbb{N}
 $$
 і  відомі  оцінки для ядер Фейєра (див., наприклад, [\ref{Z}, с. 148--151])
$$
0<F_{k}(t)<k+1, \ F_{k}(t)\leq \frac{A_{1}}{(k+1)t^{2}}, \
0<t\leq\pi,
$$
 одержимо наступні оцінки для $V_{m}(t)$:
$$
|V_{m}(t)|<A_{2}m, \ \ |V_{m}(t)|\leq\frac{A_{3}}{mt^{2}}, \ \
0<t\leq\pi,
$$
$A_{i}$ --- абсолютні сталі. Тоді, при $1\leq p<\infty$
\begin{equation}\label{na}
\|V_{2n}(t)-V_{n}(t)\|_{p}\leq
A_{4}\Bigg(\int\limits_{0\leq|t|\leq\frac{1}{n}}n^{p}dt+\int\limits_{\frac{1}{n}
\leq|t|\leq\pi}\frac{1}{(nt^{2})^{p}}dt\Bigg)^{\frac{1}{p}}\leq
A_{5}n^{1-\frac{1}{p}}.
\end{equation}

Зауважимо, що при $p=1$ з нерівності (\ref{na}) випливає оцінка
\begin{equation}\label{na1}
\|V_{2n}(t)-V_{n}(t)\|_{1}\leq A_{5}.
\end{equation}

Далі, враховуючи включення $\psi\in B$ та нерівність (\ref{na}), з
(\ref{w1}) отримуємо
\begin{equation}\label{w2}
\Big\|\big(f^{\ast}_{n,\alpha}(t)\big)^{\psi}_{\beta}\Big\|_{p}\leq\frac{\alpha
\widetilde{C}_{p,\psi}}{n^{1-\frac{1}{p}}} n^{1-\frac{1}{p}}=\alpha \widetilde{C}_{p,\psi},
\end{equation}
$\widetilde{C}_{p,\psi}$ --- величина, що залежить від $\psi$ і $p$. При
$\alpha=\alpha_{\ast}=(\widetilde{C}_{p,\psi})^{-1}$ з (\ref{w2}) випливає
нерівність (\ref{n}), а, отже, і  включення
$f^{\ast}_{n,\alpha^{\ast}}\in C^{\psi}_{\beta,p}$.

Знайдемо коефіцієнти Фур'є функції $V_{2n}(t)-V_{n}(t)$. Згідно з
формулою (3.3.5) роботи [\ref{S1}, с. 31] для ядер $V_{m}(t)$
виконується рівність
\begin{equation}\label{uad}
V_{m}(t)=D_{m}(t)+2\sum\limits_{k=m+1}^{2m-1}\Big(1-\frac{k}{2m}\Big)\cos
kt, \ m\in \mathbb{N}.
\end{equation}
Застосовуючи (\ref{uad}) при $m=n$ і $m=2n$, одержуємо
\begin{equation}\label{kfa}
V_{2n}(t)-V_{n}(t)=
$$
$$
=D_{2n}(t)-D_{n}(t)-2\sum\limits_{k=n+1}^{2n-1}\Big(1-\frac{k}{2n}\Big)\cos
kt+2\sum\limits_{k=2n+1}^{4n-1}\Big(1-\frac{k}{4n}\Big)\cos kt=
$$
$$
=-\sum\limits_{k=n+1}^{2n}\Big(1-\frac{k}{n}\Big)\cos
kt+2\sum\limits_{k=2n+1}^{4n-1}\Big(1-\frac{k}{4n}\Big)\cos kt.
\end{equation}
В силу  рівності Парсеваля
\begin{equation}\label{na2}
\|V_{2n}(t)-V_{n}(t)\|^{2}_{2}=\pi\Bigg(\sum\limits_{k=n+1}^{2n}\Big(1-\frac{k}{n}\Big)^{2}+4\sum\limits_{k=2n+1}^{4n-1}\Big(1-\frac{k}{4n}\Big)^{2}\Bigg)=
$$
$$
=\frac{\pi}{n^{2}}\Bigg(\sum\limits_{k=n+1}^{2n}(n-k)^{2}+\frac{1}{4}\sum\limits_{k=2n+1}^{4n-1}(4n-k)^{2}\Bigg)=\frac{\pi}{n^{2}}\Bigg(
\sum\limits_{k=1}^{n}k^{2}+\frac{1}{4}\sum\limits_{k=1}^{2n-1}k^{2}\Bigg)=
$$
$$
=\frac{\pi}{n^{2}}\Bigg(\frac{n(n+1)(2n+1)}{6}+\frac{(2n-1)2n(4n-1)}{24}\Bigg)=\pi(n+\frac{1}{4n}).
\end{equation}
Із  (\ref{kfa}) випливає, що $(V_{2n}-V_{n})\perp t_{n-1}$ для
будь--якого полінома
 $t_{n-1}\in~ \mathcal{T}_{2n-1}$.
 Тому
\begin{equation}\label{1bik}
\int\limits_{-\pi}^{\pi}\big(f^{\ast}_{n,\alpha_{\ast}}(t)-t_{n-1}(t)\big)\big(V_{2n}(t)-V_{n}(t)\big)dt=
\int\limits_{-\pi}^{\pi}f^{\ast}_{n,\alpha_{\ast}}(t)\big(V_{2n}(t)-V_{n}(t)\big)dt-
$$
$$
-\int\limits_{-\pi}^{\pi}t_{n-1}(t)\big(V_{2n}(t)-V_{n}(t)\big)dt=
\int\limits_{-\pi}^{\pi}f^{\ast}_{n,\alpha_{\ast}}(t)\big(V_{2n}(t)-V_{n}(t)\big)dt=
$$
$$
=\frac{\alpha_{\ast}\psi(n)}{n^{1-\frac{1}{p}}}\|V_{2n}(t)-V_{n}(t)\|^{2}_{2}.
\end{equation}
З іншого боку, використовуючи нерівність Гельдера та враховуючи
(\ref{na1}), отримуємо
\begin{equation}\label{2bik}
\int\limits_{-\pi}^{\pi}\big(f^{\ast}_{n,\alpha_{\ast}}(t)-t_{n-1}(t)\big)\big(V_{2n}(t)-V_{n}(t)\big)
dt\leq
$$
$$
\leq\|f^{\ast}_{n,\alpha_{\ast}}(t)-t_{n-1}(t)\|_{\infty}\|V_{2n}(t)-V_{n}(t)\|_{1}\leq
 A_{5}\|f^{\ast}_{n,\alpha_{\ast}}(t)-t_{n-1}(t)\|_{\infty}.
\end{equation}
Із  (\ref{na2}) -- (\ref{1bik})  одержуємо
\begin{equation}\label{ozv22}
\|f^{\ast}_{n,\alpha_{\ast}}(t)-t_{n-1}(t)\|_{\infty}\geq
\frac{\alpha_{\ast}\psi(n)}{A_{5}n^{1-\frac{1}{p}}}
\|V_{2n}(t)-V_{n}(t)\|^{2}_{2}\geq\frac{K^{(2)}_{\psi,p}\psi(n)}{n^{1-\frac{1}{p}}}
n=K^{(2)}_{\psi,p}\psi(n)n^{\frac{1}{p}}.
\end{equation}
 З  (\ref{ozv}) та
(\ref{ozv22}) випливає (\ref{t11}). Теорему 1 доведено.

\noindent \rm {\bf Теорема 2.} \it Нехай $\beta\in \mathbb{R}$,
 $\psi\in B\cap\Theta_{1}$ і
виконується одна з умов
\begin{equation}\label{umop1}
\Delta^{2}\big(1/\psi(k)\big)\geq 0, \ \ k\in\mathbb{N}
\end{equation}
або
\begin{equation}\label{umop2}
\Delta^{2}\big(1/\psi(k)\big)\leq 0, \ \ k\in\mathbb{N},
\end{equation}
де
$\Delta^{2}\big(1/\psi(k)\big)\mathop{=}\limits^{\rm df}\frac{1}{\psi(k)}-\frac{2}{\psi(k+1)}+\frac{1}{\psi(k+2)}$.

Тоді існують додатні величини $K^{(1)}_\psi$ і $K^{(2)}_\psi$, що можуть
залежати лише від $\psi$ такі, що для довільних $n\in \mathbb{N}$
\begin{equation}\label{t3}
K^{(2)}_\psi\psi(n)n\leq E_{n}\big(C^{\psi}_{\beta,1}\big)_{C}\leq{\cal
E}_{n}\big(C^{\psi}_{\beta,1})_{C}\leq K^{(1)}_\psi\psi(n)n.
\end{equation}

\bf Доведення. \ \rm Застосувавши   нерівність (\ref{ozv})
 при  $p=1$ маємо
\begin{equation}\label{ozv3}
E_{n}\big(C^{\psi}_{\beta,1}\big)_{C}\leq{\cal
E}_{n}\big(C^{\psi}_{\beta,1}\big)_{C}\leq K^{(1)}_\psi\psi(n)n,
\end{equation}
де $K^{(1)}_\psi$ --- величина, що залежить лише від $\psi$.

Для того, щоб одержати оцінку знизу, розглянемо функцію
$$
f^{(1)}_{n,\alpha}(t)=\alpha\psi(n)\big(V_{2n}(t)-V_{n}(t)\big), \
\alpha>0, \ n\in \mathbb{N},
$$
де $V_{m}(t)$ --- ядра Валле Пуссена вигляду  (\ref{uavp}).
Покажемо, що при певному виборі параметра $\alpha>0$
$f^{(1)}_{n,\alpha}\in C^{\psi}_{\beta,1}$. Для цього нам буде
корисним твердження роботи  [\ref{S2}, с. 120] в якому встановлено
оцінки норм нерівності Бернштейна для $(\psi,\beta)$-похідних в $L_1$-метриці для поліномів
$t_{m}\in\mathcal{T}_{2m+1}$.

\noindent \rm {\bf Твердження 2.} \it Нехай $\psi(k)$ --- довільна незростаюча послідовність невід'ємних чисел, для яких виконується одна з умов {\rm(\ref{umop1})} або {\rm(\ref{umop2})} і, крім того,
\begin{equation}\label{nerb2'}
\Big|\sin\frac{\beta\pi}{2}\Big|\sum\limits_{k=1}^{n-1}\psi(n)\big(k\psi(k)\big)^{-1}=O(1),
\end{equation}
де $O(1)$ --- величина, рівномірно обмежена по $n$.  Тоді для довільного тригонометричного полінома
$t_{m}(\cdot)$ порядку $m$ знайдеться стала $C_{\psi}$, що може
залежати тільки від функції від $\psi(\cdot)$, така, що
\begin{equation}\label{nerb2}
\big\|\big(t_{m}(\cdot)\big)^{\psi}_{\beta}\big\|_{1}\leq
C_{\psi}\big(\psi(m)\big)^{-1}\|t_{m}(\cdot)\|_{1}.
\end{equation}

\rm Зауважимо, що коли $\psi\in\Theta_p, \ 1\leq p<\infty$, то умова (\ref{nerb2'}) завжди виконується, оскільки в цьому випадку існує число $\alpha>\frac{1}{p}$ таке, що послідовність $\varphi(n)=n^{-\alpha}\psi(n)$ монотонно спадна, а також
\begin{equation}\label{nerb2''}
\sum\limits_{k=1}^{n-1}\frac{\psi(n)}{k\psi(k)}=\sum\limits_{k=1}^{n-1}\frac{\varphi(n)k^\alpha}{n^\alpha\varphi(k)k}\leq
\frac{K_1}{n^\alpha}\sum\limits_{k=1}^{n-1}\frac{k^\alpha}{k}\leq K_2<\infty.
\end{equation}
Зокрема, якщо $\psi\in\Theta_1$, то умова (\ref{nerb2'}) виконується при будь-яких $\beta\in\mathbb{R}$. Зогляду на це, оскільки $f^{(1)}_{n,\alpha}(t)$ --- тригонометричний поліном порядку $4n-1$, то в силу (\ref{nerb2}), з урахуванням включення
 $\psi\in B\cap \Theta_1$ та виконання однієї з умов (\ref{umop1})  або (\ref{umop2}), одержимо
\begin{equation}\label{nw2}
\Big\|\big(f^{(1)}_{n,\alpha}(t)\big)^{\psi}_{\beta}\Big\|_{1}=\alpha\psi(n)\|\big(V_{2n}(t)-V_{n}(t)\big)^{\psi}_{\beta}\|_{1}\leq
$$
$$
\leq\alpha C_{\psi}\frac{\psi(n)}{\psi(4n-1)}
\|V_{2n}(t)-V_{n}(t)\|_{1}\leq \alpha \widetilde{C}_{\psi},
\end{equation}
де $\widetilde{C}_{\psi}$ --- величина, що залежить від $\psi$. При
$\alpha=\alpha_{1}=(\widetilde{C}_{\psi})^{-1}$ з (\ref{nw2}) випливає, що
$\Big\|\big(f^{(1)}_{n,\alpha_{1}}(t)\big)^{\psi}_{\beta}\Big\|_{1}\leq
1$.

Провівши ті ж міркування, які використовувались для знаходження
оцінки (\ref{ozv22})  для функції $f^{(1)}_{n,\alpha}(t)$
одержимо нерівність
\begin{equation}\label{ozvt3}
\|f^{(1)}_{n,\alpha_{1}}(t)-t_{n-1}(t)\|_{\infty}\geq
K^{(2)}_\psi\psi(n)n,
\end{equation}
де $K^{(2)}_\psi$ --- величина, що залежить від $\psi$. Із (\ref{ozv3})
та (\ref{ozvt3}) випливає оцінка (\ref{t3}). Теорему 2 доведено.

\noindent \rm {\bf Теорема 3.} \it Нехай $1<p\leq\infty$,
$\beta\in \mathbb{R}$,
 $\psi\in B\cap\Theta_{p'}$, $\frac{1}{p}+\frac{1}{p'}=1$, і
виконується одна з умов {\rm(\ref{umop1})} або {\rm(\ref{umop2})}. Тоді
існують додатні величини $K^{(3)}_{\psi,p}$ і $K^{(4)}_{\psi,p}$, що можуть
залежати лише  від $\psi$ і $p$, такі, що для довільних $n\in
\mathbb{N}$
\begin{equation}\label{t4}
K^{(4)}_{\psi,p}\psi(n)n^{\frac{1}{p'}}\leq
E_{n}\big(L^{\psi}_{\beta,1}\big)_{p}\leq{\cal
E}_{n}\big(L^{\psi}_{\beta,1})_{p}\leq
K^{(3)}_{\psi,p}\psi(n)n^{\frac{1}{p'}}.
\end{equation}

\bf Доведення. \ \rm Зауважимо, що за виконання умови $\psi\in
\Theta_{p'}, \ 1\leq p'<\infty$, твірне ядро $\Psi_{\beta}$ класу $L^{\psi}_{\beta,1}$
задовольняє включенню $\Psi_{\beta}\in L_{p}$, $1<p\leq\infty$. Тоді $L^{\psi}_{\beta,1}\subset L_{p}$ і, з
урахуванням нерівності Юнга  (див., наприклад, [\ref{S1}, с. 293]) та
інтегрального зображення (\ref{zgo}), маємо
\begin{equation}\label{intzob4}
{\cal E}_{n}\big(L^{\psi}_{\beta,1}\big)_{p}\leq
\frac{1}{\pi}\big\|\Psi_{\beta,n}(\cdot)\big\|_{p}\|
\varphi(\cdot)\|_{1}\leq\frac{1}{\pi}\big\|\Psi_{\beta,n}(\cdot)\big\|_{p}.
\end{equation}

В силу (\ref{ny})
\begin{equation}\label{np'}
\big\|\Psi_{\beta,n}(\cdot)\big\|_{p}\leq
K_{p',1}\big(\sum\limits_{k=n}^{\infty}\Delta\psi(k)k^{\frac{1}{p'}}+\psi(n)n^{\frac{1}{p'}}\big),
\ \ 1< p\leq\infty, \ \frac{1}{p}+\frac{1}{p'}=1,
\end{equation}
$K_{p,1}$ --- стала, що залежить від $p$.
 Тоді, застосувавши до суми $\sum\limits_{k=n}^{\infty}\Delta\psi(k)k^{\frac{1}{p'}}$  лему 1,
із (\ref{intzob4}) та (\ref{ny}) одержимо
\begin{equation}\label{t4ozver}
E_{n}\big(L^{\psi}_{\beta,1}\big)_{p}\leq{\cal
E}_{n}\big(L^{\psi}_{\beta,1}\big)_{p}\leq
K^{(3)}_{\psi,p}\psi(n)n^{\frac{1}{p'}}.
\end{equation}

Щоб одержати оцінку знизу величини
$E_{n}\big(L^{\psi}_{\beta,1}\big)_{p}$, $1<p\leq\infty$,
розглянемо функцію
$$
f^{(1)}_{n,\alpha}(t)=\alpha\psi(n)\big(V_{2n}(t)-V_{n}(t)\big), \
\alpha>0, \ n\in \mathbb{N}.
$$
Як випливає з нерівності (\ref{nw2}), при певному виборі параметра
$\alpha=\alpha_{1}$, залежного від $\psi$, виконуватиметься
нерівність
$\Big\|\big(f^{(1)}_{n,\alpha_{1}}(t)\big)^{\psi}_{\beta}\Big\|_{1}\leq
1$ і, отже, $f^{(1)}_{n,\alpha_{1}}\in L^{\psi}_{\beta,1}$.

Оскільки $(V_{2n}-V_{n})\perp t_{n-1}$ для будь-якого полінома $t_{n-1}\in\mathcal{T}_{2n-1}$, то   має місце рівність
\begin{equation}\label{ocht41}
\int\limits_{-\pi}^{\pi}\big(f^{(1)}_{n,\alpha_{1}}(t)-t_{n-1}(t)\big)\big(V_{2n}(t)-V_{n}(t)\big)dt=
\alpha_{1}\psi(n)\|V_{2n}(t)-V_{n}(t)\|^{2}_{2}.
\end{equation}
З іншого боку, в силу  (\ref{na})
$$
\|V_{2n}(t)-V_{n}(t)\|_{p'}\leq A_{5}n^{\frac{1}{p}}, \ \
1<p\leq\infty,
$$
і тому, застосовуючи нерівність 3.8.1 та 3.8.3 роботи [\ref{S1}, с. 137, 138],
маємо

\begin{equation}\label{ocht42}
\int\limits_{-\pi}^{\pi}\big(f^{(1)}_{n,\alpha_{1}}(t)-t_{n-1}(t)\big)\big(V_{2n}(t)-V_{n}(t)\big)dt\leq
\|f^{(1)}_{n,\alpha_{1}}(t)-t_{n-1}(t)\|_{p}
\|V_{2n}(t)-V_{n}(t)\|_{p'}\leq
$$
$$
\leq
A_{5}n^{\frac{1}{p}}\|f^{(1)}_{n,\alpha_{1}}(t)-t_{n-1}(t)\|_{p}.
\end{equation}
Із (\ref{ocht41}) та (\ref{ocht42}), враховуючи співвідношення
(\ref{na2}), отримуємо
\begin{equation}\label{ochf4}
\|f^{(1)}_{n,\alpha_{1}}(t)-t_{n-1}(t)\|_{p}\geq\frac{\alpha_{1}\psi(n)}{A_{5}n^{\frac{1}{p}}}\|V_{2n}(t)-V_{n}(t)\|^{2}_{2}=
\frac{K^{(4)}_{\psi,p}\psi(n)}{n^{\frac{1}{p}}}n=K^{(4)}_{\psi,p}\psi(n)n^{\frac{1}{p'}}.
\end{equation}
З (\ref{t4ozver}) та (\ref{ochf4}) випливає (\ref{t4}). Теорему 3
доведено.

Зауважимо, що в теоремах 2 і 3 вимоги виконання однієї з умов {\rm(\ref{umop1})} та {\rm(\ref{umop2})} можна замінити на більш загальну (але менш прозору): щоб для $\beta\in\mathbb{R}$ і для послідовності ${\psi\in B\cap\Theta_1}$ (у випадку теореми 2) або $\psi\in B\cap\Theta_p$ (у випадку теореми 3) виконувалась нерівність (\ref{nerb2}).

 Для функцій $\psi(t)=\frac{\ln^{\alpha}(t+c)}{t^{r}}$,
$\alpha\geq1$,
 $c>e^{\frac{2\alpha}{r}}-1$, $t\geq 1$ та
$\psi(t)=\frac{1}{t^{r}\ln^{\alpha}(t+c)}$, $\alpha\geq 0$, $c>0$,
$t\geq1$ виконуються всі умови теорем 2 (при $r>\frac{1}{p}$) та 3
(при $r>1-\frac{1}{p}$) і, отже, для величин $E_n(C^\psi_{\beta,1})_C$ та $E_n(L^\psi_{\beta,1})_p, \ 1<p\leq\infty$, мають місце співвідношення (\ref{t3}) та (\ref{t4}) відповідно.
\newpage

\begin{enumerate}
\Rus

\centerline {\bf Література}

\item \label{S1}
{\it Степанец А.И.} Методы теории приближений. --- Киев: Ин-т
математики НАН Украины, 2002. --- {\bf 40}. --- Ч.І. --- 427 с.

\item\label{Korn}
{\it Корнейчук Н.П.} Экстремальные задачи теории приближения.
--- М.: Наука, 1976.
--- 320 с.

\item\label{T}
{\it Temlyakov V.N.} Approximation of Periodic Function: Nova
Science Publi--\ chers, Inc. --- 1993. --- 419p.

\item \label{Kol}
{\it Kolmogoroff A.} Zur Gr\"{o}ssennordnung des Restgliedes
Fourierschen Reihen differenzierbarer Funktionen // Ann. Math.(2),
--- 1935. --- {\bf 36}, №2. --- P.~ 521--526.

\item \label{Pin}
{\it Пинкевич В.Т.} О порядке остаточного члена ряда Фурье
функций, дифференцируемых в смысле Вейля // Изв. АН СССР. Сер.
мат. --- 1940. --- {\bf 4}, №6. --- С. 521--528.

\item\label{Nik}
{\it Никольский С.М.} Приближение периодических функций
тригонометрическими многочленами // Труды МИАН СССР. --- 1945. ---
{\bf 15}. --- С. 1--76.

\item\label{Nik2}
{\it Никольский С.М.} Приближение  функций тригонометрическими
полиномами в среднем
// Изв. АН
СССР. Cер. мат.
--- 1946. --- {\bf 10}, №3. --- С. 207--256.

\item \label{Ef}
{\it Ефимов А.В.} Приближение непрерывных периодических функций
суммами Фурье // Изв. АН СССР Сер. мат. --- 1960. --- {\bf 24},
№2.
--- С. 243--296.

\item \label{Tel}
{\it Теляковский С.А.} О нормах тригонометрических полиномов и
приближении дифференцируемых функций линейными средними их рядов
Фурье // Труды МИАН СССР. --- 1961. --- {\bf 62}. --- С. 61--97.

\item\label{Fav}
{\it Favard J.} Sur l'approximation des fonctions p\'{e}riodiques
par des polynomes trigonom\'{e}triques
// C.R. Acad. Sci. --- 1936. --- {\bf 203}. --- P. 1122--1124.

\item\label{Fav1}
{\it Favard J.} Sur les meilleurs proc\'{e}des d'approximations de
certains classes de  fontions par des polynomes
 trigonom\'{e}triques // Bull. de Sci. Math. ---
1937. --- {\bf 61}. --- P. 209--224, 243--256.

\item\label{Dz}
{\it Дзядык В.К.} О наилучшем приближении на классе периодических
функций, имеющих ограниченную $s$--ю производную $(0<s<1)$ // Изв.
АН СССР, Сер. мат. --- 1953. --- {\bf 17}. --- С. 135--162.

\item\label{Dz1}
{\it Дзядык В.К.} О наилучшем приближении на классах периодических
функций, определяемых интегралами от линейной комбинации абсолютно
монотонных ядер  // Мат. заметки. --- 1974. --- {\bf 16}, №5.
--- С.~ 691--701.

\item\label{ST}
{\it Стечкин С.Б.} О наилучшем приближении некоторых классов
периодических функций тригонометрическими полиномами // Изв. АН
СССР, Cер. мат. --- 1956. --- {\bf 20}, --- С. 643--648.

\item\label{SU}
{\it Сунь Юн--шен} О наилучшем приближении периодических
дифференцируемых функций тригонометрическими полиномами // Изв. АН
СССР. Cер. мат.
--- 1959. --- {\bf 23}, №1. --- С. 67--92.

\item\label{Bab}
{\it Бабенко В.Ф., Пичугов С.А.}  О наилучшем линейном
приближении некоторых классов дифференцируемых периодических
 функций // Мат. заметки. --- 1980. --- {\bf 27}, №5. --- С.~ 683--689.

\item\label{Kuchpel}
{\it Степанец А.И., Кушпель А.К.}  Скорость сходимости рядов Фурье
и наилучшие приближения в пространстве $L_{p}$
// Укр. мат. журн. --- 1987.
--- {\bf 39}, №4. --- С. 483--492.

\item \label{Step monog 1987} {\it  Степанец А.И.} Классификация и приближение периодических функций. --- Киев: Наук. думка~--- 1987.~--- 268 c.

\item \label{Rom}
{\it Романюк В.С.} Дополнения к оценкам приближения суммами Фурье
классов бесконечно дифференцируемых функций // Екстремальні задачі
теорії функцій та суміжні питання: Праці Ін-ту математики НАН
України.
--- 2003. --- {\bf 46},
 --- С. 131--135.

\item \label{Serduk}
{\it Сердюк А.С., Соколенко І.В.} Рівномірні наближення класів
$(\psi, \overline{\beta})$--диференційовних функцій лінійними
методами
// Зб. праць Ін-ту матем. НАН України. --- 2011. --- {\bf 8}, №1.
 --- С. 181--189.

 \item \label{S2}
{\it Степанец А.И.} Методы теории приближений. --- Киев: Ин-т
математики НАН Украины, 2002. --- {\bf 40}. --- Ч.II. --- 468 с.

\item \label{Bari}
{\it Бари Н.К.} Тригонометрические ряды. --- М.: Физматгиз, 1961.
 --- 936~с.

\item\label{Z}
{\it Зигмунд А.} Тригонометрические ряды. В 2 т.~--- М.: Мир,
1965.~--- Т.І.~ --- 615 с.


\end{enumerate}
\end{document}